\documentclass{ifacconf}

\usepackage{graphicx}      
\usepackage{natbib}        

\usepackage{amsmath,amssymb}
\begin{document}
\begin{frontmatter}

\title{Obstruction to the billinear control of the Gross-Pitaevskii equation: an example with an unbounded potential\thanksref{footnoteinfo}} 

\thanks[footnoteinfo]{T. Chambrion is  supported by the grant "QUACO" ANR-17-CE40-0007-01. L. Thomann is supported by the grants "BEKAM"  ANR-15-CE40-0001 and "ISDEEC'' ANR-16-CE40-0013.}

\author[First]{Thomas Chambrion} 
\author[Second]{Laurent Thomann}

\address[First]{Universit\'e de Lorraine, CNRS, Inria, IECL, F-54000 Nancy, France (e-mail: Thomas.Chambrion@univ-lorraine.fr).}
\address[Second]{Universit\'e de Lorraine, CNRS, IECL, F-54000 Nancy, France (e-mail: Laurent.Thomann@univ-lorraine.fr). }

\begin{abstract}                
In 1982, Ball, Marsden, and Slemrod proved an obstruction to the controllability 
of linear dynamics with a bounded bilinear control term.  This note presents an 
example of nonlinear dynamics  with respect to the state for which this 
obstruction still holds while the control potential is not bounded. 
\end{abstract}


\end{frontmatter}

\section{Introduction and results} 

\subsection{Introduction}

On the Euclidean space $\mathbf{R}^3$ endowed with its natural norm $| \cdot |$, we study the   control  problem: 
\begin{equation}\label{NLS}
  \left\{
      \begin{array}{ll}
      & i \partial_t \psi +H \psi =u(t) K(x) \psi- \sigma |\psi|^2 \psi, \! \!\!\!\!\qquad (t,x)\in \mathbf{R} \times \mathbf{R}^3,
       \\  & \psi(0,x)=\psi_0(x),
      \end{array}
    \right.
\end{equation} 
where  $H=-\Delta+ |x|^2$ is the Hamiltonian of the quantum harmonic oscillator on $\mathbf{R}^3$, $u : \mathbf{R} \longrightarrow \mathbf{R}$ is the control, $K: \mathbf{R}^3 \longrightarrow \mathbf{R}$ is a given potential and $\sigma \in \{0,1\}$. 

The  Sobolev spaces  based on the domain of the harmonic oscillator are instrumental in the study of dynamics (\ref{NLS}). They are defined, for $s\geq 0$ and $p\geq 1$ by
 \begin{equation*} 
         {\mathcal{W}}^{s, p}= {\mathcal{W}}^{s, p}(\mathbf{R}^3) = \big\{ f\in L^p(\mathbf{R}^3),\; {H}^{s/2}f\in L^p(\mathbf{R}^3)\big\},     
       \end{equation*}
       \begin{equation*}
           {\mathcal{H}}^{s}=   {\cal H}^{s}(\mathbf{R}^3) = {\mathcal{W}}^{s, 2}.
       \end{equation*}
             The natural norms are denoted by $\Vert f\Vert_{{\mathcal{W}}^{s,p}}$ and up to equivalence of norms  (see {\it e.g.} Lemma~2.4 of \cite{YajimaZhang2}), for $1<p<+\infty$,  we have
          \begin{equation}\label{equiv}
      \Vert f\Vert_{{\mathcal{W}}^{s,p}} = \Vert  H^{s/2}f\Vert_{L^{p}} \equiv \Vert (-\Delta)^{s/2} f\Vert_{L^{p}} + 
       \Vert\langle x\rangle^{s}f\Vert_{L^{p}},
 \end{equation}
 with the notation $\langle x\rangle=(1+|x|^{2})^{1/2}$.\medskip

\subsection{Ball-Marsden-Slemrod obstructions}

The dynamical system \eqref{NLS} is called the bilinear Gross-Pitaevskii equation. 
It is a nonlinear version of 
dynamics of the type $\dot{\psi}=A\psi + u(t) B\psi$ where $A$ and $B$ are linear 
operators in a Banach space $X$ and $u:\mathbf{R} \to \mathbf{R}$ is a real-valued control  which involve a  control term $u B \psi$ that is bilinear in 
$(u,\psi)$. Such dynamics play a major role in physics and are the subject of a
vast literature \cite{Khapalov}.   In \cite{BMS}, Ball, Marsden, and Slemrod have proven that if $A$ 
generates a $C^0$ semi-group in $X$ and if $B$ is bounded on $X$, then the 
attainable set from any source $\psi_0$ in $X$ with $L^r$ controls, $r>1$, is 
contained in a countable union of compact sets of~$X$. This represents a deep obstruction to the controllability of bilinear control systems in infinite dimensional Banach spaces, since this result implies that the attainable set is meager in Baire sense and has empty interior.

The original result of \cite{BMS} (and its adaptation to the Schr\"odinger 
equation in \cite{Turinici}) has been extended to the case of $L^1$ controls in 
\cite{BCC2}. More recently, the case where $A$ is non-linear has been investigated in \cite{ChTh1} (for the Klein-Gordon 
equation)  and \cite{ChTh2} (for the Gross-Pitaevskii equation~\eqref{NLS}).

%

In \cite[Theorem 1.6]{ChTh2} we showed in particular that if $K \in {\mathcal{W}}^{1,\infty}(\mathbf{R}^3)$, the 
dynamics~\eqref{NLS} is non controllable. Under this assumption,  the map 
 \begin{equation*} 
 \begin{array}{rcl}
{\mathcal{H}}^1(\mathbf{R}^3)&\longrightarrow&{\mathcal{H}}^1(\mathbf{R}^3)\\[3pt]
\displaystyle  \psi&\longmapsto & K\psi,
 \end{array}
 \end{equation*}
 is continuous and this was used in the heart of the  proof. 
 
 The main result of this note, Theorem~\ref{BMS_NLS2} below, provides an example of potential $K \notin L^{\infty}(\mathbf{R}^3)$ where this condition is violated, but where the  obstruction to controllability result still holds true. 

\subsection{Main result}

Our main result reads as follows

     \begin{thm}\label{BMS_NLS2}
Let   $K(x)=\log (|x|){\bf 1}_{\{|x|\leq 1\}}$ and  $\psi_0 \in {\mathcal{H}}^1(\mathbf{R}^3) $. Assume that $u \in \bigcup_{r>1}L^r_{loc}(\mathbf{R})$, then the equation~\eqref{NLS} admits a global flow  $\psi(t)=\Phi^u(t)(\psi_0)$. 

Moreover, for every $\psi_0 \in {\mathcal{H}}^1(\mathbf{R}^3) $, the  attainable set 
\begin{equation}\label{EQ_attainable_set}
\bigcup_{t \in \mathbf{R}} \;\bigcup_{\substack{ u \in L^r_{loc}(\mathbf{R}), \\  r>1}} \big\{\Phi^u(t) (\psi_0)\big\}
\end{equation}
is   a countable union of compact subsets of ${\mathcal{H}}^1(\mathbf{R}^3)$.
\end{thm}
  In this paper, the solutions to~\eqref{NLS} are understood in the mild sense
    \begin{eqnarray*} 
 \lefteqn{\psi(t)=e^{it H}\psi_0-i \int_0^t u(\tau) e^{i(t-\tau) H}( K \psi(\tau))d\tau}\\
 & ~~~~~~~~~~~~~~~~~& +i \sigma \int_0^t e^{i(t-\tau) H}(|\psi|^2 \psi)d\tau.
 \end{eqnarray*} 

%

 \subsection{Content of the paper} The rest of this note provides a  proof of Theorem \ref{BMS_NLS2}.  The proof crucially relies on classical Strichartz estimates, which we recall in  Section  \ref{SEC_Strichartz}. The proof itself is split in two parts. The global well-posedness of the problem \eqref{NLS} is established in Section \ref{SEC_Well_posedness}, using among other some energy estimates. The proof of the obstruction result follows the strategy used in the paper \cite{BMS} and is given in Section \ref{SEC_meagerness}.

\section{Strichartz estimates} \label{SEC_Strichartz}
 
 As in \cite{ChTh2}, the  Strichartz estimates play a major role in the argument, let us recall them in the three-dimensional case.
A couple $(q,r)\in [2,+\infty]^2$ is called admissible if 
\begin{equation*}
\frac2q+\frac{3}{r}=\frac{3}2,
\end{equation*}
and if one defines 
\begin{equation*} 
 {X}^1_{T}:= \bigcap_{(q,r) \ admissible} L ^q\big( [-T, T] \,; {\mathcal{W}}^{1,r}( \mathbf{R}^3 )\big),
\end{equation*}
then for all $T>0$ there exists $C_{T}>0$ so that for all $\psi_{0}\in {\mathcal{H}}^{1}(\mathbf{R}^3)$ we have 
\begin{equation}\label{Stri}
\|\e^{itH}\psi_{0}\|_{X^{1}_{T}}\leq C_{T}\|\psi_{0}\|_{{\mathcal{H}}^{1}(\mathbf{R}^3).}
\end{equation}
Using interpolation theory one can prove that 
\begin{equation*} 
 {X}^1_{T}=   L ^{\infty}\big( [-T, T] \,; {\mathcal{H}}^1( \mathbf{R}^3 )\big) \cap  L ^{2}\big( [-T, T] \,; {\mathcal{W}}^{1,6}( \mathbf{R}^3 )\big),
\end{equation*}
so that one can define
$$\|\psi\|_{X^1_T}= \| \psi \|_{L^{\infty}([-T,T];{\mathcal{H}}^1(\mathbf{R}^3))}+\| \psi \|_{L^{2}([-T,T];{\mathcal{W}}^{1,6}(\mathbf{R}^3))}.$$

We will also need the inhomogeneous version of Strichartz:  for all $ T>0 $, there exists $C_{T}>0$ so that for all   admissible couple       $ ( q, r ) $ and function  $ F \in L^{q'}( [T,T]; {\mathcal{W}}^{1,r'} (\mathbf{R}^3)) $,
\begin{equation}\label{Stri1}
 \big\|  \int _0^t \e^{i(t-\tau)H} F(\tau) d\tau   \big\| _{   X^1_T} \leq C_{T} \| F \|_{  L^{q'} ([-T,T],{{\mathcal{W}}}^{1,r'} (\mathbf{R}^3)) },
\end{equation}
where $ q' $ and $ r'$ are the H\"older conjugate of $ q $ and $ r $. We refer to \cite{poiret2} for a proof. Let us point out that \eqref{Stri1} implies that
\begin{eqnarray}\label{Stri1.1}
 \lefteqn{\big\|  \int _0^t \!\!\! \e^{i(t-\tau)H} F(\tau) d\tau   \big\| _{   X^1_T} \leq}
 \nonumber \\
 && C_{T} \big(\| F_1 \|_{  L^{1} ([-T,\!T],{{\mathcal{H}}}^{1} (\mathbf{R}^3)) }\!+\!\| F_2 \|_{  L^{2} ([-T,\!T],{{\mathcal{W}}}^{1,\frac65} (\mathbf{R}^3)) } \big)
\end{eqnarray}
for any $F_1, F_2$ such that $F_1+F_2=F$, which will prove useful.
\medskip

In the sequel $c,C>0$ denote constants the value of which may change
from line to line. These constants will always be universal, or uniformly bounded. For $x \in \mathbf{R}^3$, we write $\langle x\rangle=(1+|x|^{2})^{1/2}$. We will sometimes use the notations $L^{p}_{T}=L^{p}([0,T])$ and $L^{p}_{T}X=L^{p}([0,T]; X)$  for $T>0$.

\section{Proof of Theorem \ref{NLS}}

             \subsection{Global existence theory for dynamics \eqref{NLS}} \label{SEC_Well_posedness}

Using the reversibility of the equation~\eqref{NLS}, it is enough to consider non-negative times in the proofs.\medskip

  The    following result will be useful to control the bilinear term in~\eqref{NLS}.

    \begin{lem}\label{lem-est}
    Let $K(x)=\log (|x|){\bf 1}_{\{|x|\leq 1\}}$ and $T>0$. Then for all $2 \leq q <\infty$,  there exists $C_T>0$ such that for all $\psi \in X^1_T$
     \begin{equation*} 
      \|K \psi \|_{L^{q}([-T,T]; {\mathcal{H}}^1(\mathbf{R}^3))} \leq C_T \|\psi \|_{X^1_T}.  
       \end{equation*}
    \end{lem}
    
    \begin{pf}
    Firstly by \eqref{equiv} we have
\begin{equation}\label{decomp}
 \big\| K \psi \big\|_{{\mathcal{H}}^1} \leq c  \big\| \nabla K \psi \big\|_{L^2} +c  \big\| K \nabla  \psi \big\|_{L^2}+c  \big\| K \langle x\rangle\psi \big\|_{L^2} .
\end{equation}

 $\bullet$ Let us study the first term in \eqref{decomp}. Since $|\nabla K | \leq C |x|^{-1}$, we can use the Hardy inequality 
            \begin{equation}\label{hardy}
           \big\|    |x|^{-1}   \psi    \big\|_{L^2(\mathbf{R}^3)} \leq C  \big\|      \psi    \big\|_{H^1(\mathbf{R}^3)}
           \end{equation} 
           (we refer to \cite[Lemma A.2]{Tao} for the general statement and proof of this inequality), and therefore the contribution of the first term reads
            $\big\| \nabla K \psi \big\|_{L^{q}_TL^2} \leq C T^{1/q} \|\psi \|_{L_T^{\infty}{\mathcal{H}}^1} \leq C_T \|\psi \|_{X^1_T}$.
            
             $\bullet$ To bound the contribution of the two last  terms in~\eqref{decomp}, we will use that ${K \in L^p(\mathbf{R}^3)}$ for any $1\leq p<\infty$. Given $2 \leq q <\infty$, we choose $r>2$   such that the couple $(q,r)$ is (Strichartz) admissible and write, using H\"older  
              $$\big\|  K \nabla\psi \big\|_{ L^2} \leq    \|  K   \|_{L^{p}}   \|  \nabla \psi  \|_{ L^r} \leq  c   \|  K   \|_{L^{p}}   \|  \psi  \|_{ {\mathcal{W}}^{1,r}},$$
    with $1/p+1/r=1/2$.         Thus 
  $$\big\|  K \nabla\psi \big\|_{ L_T^{q} L^2} \leq    c   \|  K   \|_{L^{p}}   \|  \psi  \|_{ L_T^{q}{\mathcal{W}}^{1,r}}\leq c   \|  K   \|_{L^{p}}    \|\psi \|_{X^1_T}.$$
             Similarly,\\  $\big\|  K \langle x\rangle\psi \big\|_{L^{q}_TL^2} \leq    c\|  K  \|_{L^{p}}   \|  \psi \|_{L^{q}_T{\mathcal{W}}^{1,r}}  \leq    c\|  K  \|_{L^{p}}     \|\psi \|_{X^1_T}$.    
    \end{pf}
    
    We now state a global existence result for \eqref{NLS} adapted to our control problem. 
    
     \begin{prop}\label{propK-nb}  Let $u\in L^{1}_{loc}(\mathbf{R})$ and set $K:x \mapsto\log (|x|){\bf 1}_{\{|x|\leq 1\}}$.
 Let $\psi_0 \in {\mathcal{H}}^1(\mathbf{R}^3)$, then the  equation~\eqref{NLS} admits a unique global solution   ${\psi \in \mathcal{C}(\mathbf{R} ; {\mathcal{H}}^1(\mathbf{R}^3))} \cap L^2_{loc}(\mathbf{R} ; {\mathcal{W}}^{1,6})$ which moreover satisfies the bounds
 \begin{eqnarray}
 \lefteqn{\| \psi \|_{L^{\infty}([-T,T];{\mathcal{H}}^1(\mathbf{R}^3))} } \nonumber \\
& ~~~~~~~~~~~~~~~\leq & C(  \| \psi_0 \|_{ {\mathcal{H}}^1(\mathbf{R}^3) }) \exp\Big(c \int_{-T}^T |u(\tau)|d\tau\Big), \label{eqB211}
\end{eqnarray}  
and
          \begin{equation} \label{boundW16}
 \|  \psi \|_{L^2([-T,T];{\mathcal{W}}^{1,6}(\mathbf{R}^3))}  \leq  C\big(  T,  \|\psi_0\|_{{\mathcal{H}}^1(\mathbf{R}^3)} ,    \int_{-T}^T |u(\tau)|   d\tau\big).
   \end{equation} 
 \end{prop}

\begin{pf}
The proof is in the spirit of the proof of \cite[Proposition 1.5]{ChTh2}, but here we use moreover the Hardy inequality~\eqref{hardy} to control the bilinear term in~\eqref{NLS}.\medskip

  {\it Energy bound:} Assume for a moment that the solution exists on a time interval $[0,T]$. For $0\leq t\leq T$, we define 
   \begin{eqnarray*} 
       E(t)&=&\int_{\mathbf{R}^3}\big(\overline{\psi} H\psi+|\psi|^2 +\frac{\sigma}2|\psi|^4\big)dx \\
       &=& \int_{\mathbf{R}^3}\big( |\nabla \psi|^2 +|x|^2|\psi|^2+|\psi|^2 +\frac{\sigma}2|\psi|^4\big)dx.
       \end{eqnarray*}
 Then, using that $  \partial_t \overline{\psi} =-i( H \overline{ \psi} +\sigma |\psi|^2 \overline{ \psi})+i u(t) K(x) \overline{ \psi}$, we get 
    \begin{eqnarray*} 
       E'(t) &=   & 2 \Re \int_{\mathbf{R}^3} \partial_t\overline{ \psi}\big(\psi+ H \psi + \sigma |\psi|^2 \psi  \big)dx \nonumber \\
       &=& -2 u(t) \Im \int_{\mathbf{R}^3}  K \overline{ \psi }H \psi dx \nonumber \\
              &=& 2 u(t) \Im \int_{\mathbf{R}^3}  \overline{ \psi }  \nabla K \cdot \nabla  \psi    dx.
                  \end{eqnarray*} 
        Observing that  $|\nabla K | \leq C |x|^{-1}$, by the Hardy inequality~\eqref{hardy} we get
       \begin{eqnarray*} 
       E'(t) &\leq&  2   |u(t)| \big\|      \psi   \nabla K \big\|_{L^2} \big\|  \nabla  \psi    \big\|_{L^2} \\
       & \leq & C   |u(t)| \big\|    |x|^{-1}   \psi    \big\|_{L^2} \big\|  \nabla  \psi    \big\|_{L^2}\\
       &\leq & C   |u(t)| \big\|     \psi    \big\|^2_{{\mathcal{H}}^1}  \\
         &\leq & C  |u(t)| E(t).
                \end{eqnarray*}        
 Thus, using that $\sigma \geq 0$, we deduce the  bound  \eqref{eqB211}.      
  
  \medskip
  {\it Local existence and global existence:} We consider the map 
  \begin{eqnarray*} 
 \lefteqn{\Phi(\psi)(t)=e^{it H}\psi_0+i \sigma \int_0^t e^{i(t-\tau) H}(|\psi|^2 \psi)d\tau}\\
 &~~~~~~~~~~~~~&-i \int_0^t u(\tau) e^{i(t-\tau) H}( K \psi)d\tau,
 \end{eqnarray*} 
and we will show that it is a contraction in the space 
  \begin{equation*}
B_{T,R}:= \big\{ \|\psi\|_{X^1_T}\leq R\big\},
 \end{equation*} 
 with $R>0$ and $T>0$ to be fixed.

 By the Strichartz inequalities \eqref{Stri} and \eqref{Stri1.1}
     \begin{eqnarray*}
     \lefteqn{\| \Phi(\psi)\|_{X^1_T} }\\
     &&  \leq c \|\psi_0\|_{{\mathcal{H}}^1}+c 
      \big\| |\psi|^2 \psi \big\|_{L^1_T{\mathcal{H}}^1} +c   \|  u   K\psi \|_{L^2_t {\mathcal{W}}^{1,6/5}(\mathbf{R}^3)}. 
     \end{eqnarray*}
Then by \cite[Lemma A.1]{ChTh2}
  \begin{eqnarray*}
\lefteqn{\|  K\psi \|_{{\mathcal{W}}^{1,6/5}(\mathbf{R}^3)} \leq C \|  K \|_{L^{3}(\mathbf{R}^3)} \|  \psi \|_{{\mathcal{H}}^{1}(\mathbf{R}^3)}}\\
&& +C \|  K \|_{{\mathcal{W}}^{1,3/2}(\mathbf{R}^3)} \|  \psi \|_{L^{6}(\mathbf{R}^3)} .
 \end{eqnarray*}
  Next, by the Sobolev embedding ${\mathcal{H}}^1(\mathbf{R}^3) \subset L^6(\mathbf{R}^3)$, from the previous line we get 
  \begin{equation}\label{borneKW}
\|  K\psi \|_{{\mathcal{W}}^{1,6/5}(\mathbf{R}^3)}  \leq C  \|  \psi \|_{{\mathcal{H}}^{1}(\mathbf{R}^3)},
 \end{equation}
 which in turn implies 
   \begin{eqnarray*} 
 \lefteqn{\| \Phi(\psi)\|_{X^1_T} }\\
 &\leq &c \|\psi_0\|_{{\mathcal{H}}^1}+c \| \psi \|_{L_T^{\infty}{\mathcal{H}}^1}\| \psi \|^2_{L^{2}_T L^{\infty}} +c\| u  \|_{L^{1}_T}  \|  \psi \|_{L^{\infty}_T{\mathcal{H}}^1}.
 \end{eqnarray*} 
 By the Gagliardo-Nirenberg and Sobolev inequalities on~$\mathbf{R}^3$, 
$$   \|  \psi\|_{ L^{\infty}} \leq C \|  \psi\|^{1/2}_{  L^6}      \|  \psi\|^{1/2}_{ {\mathcal{W}}^{1,6}} \leq C     \|  \psi\|^{1/2}_{  {\mathcal{H}}^1}      \|  \psi\|^{1/2}_{ {\mathcal{W}}^{1,6}}   , $$
 thus    $  \| \psi \|_{L^{2}_T L^{\infty}}  \leq cT^{1/4}\|  \psi \|^{1/2}_{L^{\infty}_T{\mathcal{H}}^1}   \|  \psi\|^{1/2}_{ L^2_T{\mathcal{W}}^{1,6}}   $ , and for $\psi \in B_{T,R}$ we get     
   \begin{equation*}
 \| \Phi(\psi)\|_{X^1_T}
  \leq c \|\psi_0\|_{{\mathcal{H}}^1}+cT^{1/2} R^3 +cR  \| u   \|_{L^{1}_T} .
 \end{equation*} 
   We now choose $R =2c \|\psi_0\|_{{\mathcal{H}}^1}$. Then for $T>0$ small enough, $\Phi$ maps $B_{T,R}$ into itself. With similar estimates we can show that $\Phi$ is a contraction in $B_{T,R}$, namely 
    \begin{equation*}
 \| \Phi(\psi_1)- \Phi(\psi_2)\|_{X^1_T}  
 \leq  \big[cT^{1/2} R^2 + c\| u   \|_{L^{1}_T}\big]  \| \psi_1- \psi_2\|_{X^1_T}.   
  \end{equation*} 
 As a conclusion there exists a unique fixed point to $\Phi$, which is a local solution to \eqref{NLS}. \medskip
 
The local time of existence only depends on $u$ and on the ${\mathcal{H}}^1$-norm. Therefore one can use the energy bound to show the global existence. \medskip

  {\it Proof of the bound \eqref{boundW16}:} The proof follows the main lines of \cite[Bound (1.18)]{ChTh2}, hence we do not detail it here.
  \end{pf}

             \subsection{Meagerness of the attainable set}\label{SEC_meagerness}

 Let $\epsilon >0$ and let $u, u_n \in L^{1+\epsilon}([0,T]; \mathbf{R})$ such 
 that $u_n \rightharpoonup u$ weakly  in $L^{1+\epsilon}([0,T]; \mathbf{R})$. 
 This implies a bound $\|u_n \|_{L^{1+\epsilon}_T} \leq C(T)$ for some $C(T)>0$, 
 uniformly in $n \geq 1$. We have
  \begin{eqnarray*} 
\lefteqn{ \psi(t)=e^{it H}\psi_0-i \int_0^t u(\tau) e^{i(t-\tau) H}( K \psi(\tau))d\tau} \\
&~~~~~~~~~~~~~&+i \sigma \int_0^t e^{i(t-\tau) H}(|\psi|^2 \psi)d\tau,
 \end{eqnarray*} 
and 
  \begin{eqnarray*} 
\lefteqn{ \psi_n(t)=e^{it H}\psi_0-i \int_0^t u_n(\tau) e^{i(t-\tau) H}( K \psi_n(\tau))d\tau}\\
&~~~~~~~~~~~~~&+i \sigma \int_0^t e^{i(t-\tau) H}(|\psi_n|^2 \psi_n)d\tau.
 \end{eqnarray*} 
We set $z_n=\psi-\psi_n$, then $z_n$ satisfies 
\begin{equation}\label{zn}
z_n= \mathcal{L}(\psi,\psi_n)+\mathcal{N}(\psi,\psi_n),
\end{equation}
with 
\begin{eqnarray*}
\lefteqn{\mathcal{L}(\psi,\psi_n)= -i   \int_0^t \big(u(\tau)-u_n(\tau) \big) e^{i(t-\tau) H}( K \psi)d\tau  }\\
&~~~~~~~~~~&-i \int_0^t  u_n(\tau)  e^{i(t-\tau) H}\big( K (\psi-\psi_n)\big)d\tau
\end{eqnarray*}
and 
\begin{eqnarray*}
\lefteqn{\mathcal{N}(\psi,\psi_n)= i \sigma \int_0^t   e^{i(t-\tau) H}\big(  (\psi-\psi_n)(\psi+\psi_n)  \overline{\psi} \big)d\tau }\\
&~~~~~~~~~~&+i \sigma \int_0^t   e^{i(t-\tau) H}\big(  (\overline{\psi}-\overline{\psi_n})\psi^2_n  \big)d\tau. 
\end{eqnarray*}
Let us prove that  $z_n \longrightarrow 0$ in $L^{\infty}([0,T]; {\mathcal{H}}^1(\mathbf{R}^3))$.  

\begin{lem}
Denote by 
  \begin{equation*} 
 \epsilon_n:=\Big\|    \int_0^t \big(u_n(\tau) -u(\tau) \big)e^{i(t-\tau) H}( K \psi(\tau))d\tau \Big\|_{L^{\infty}_T{\mathcal{H}}^1(\mathbf{R}^3)}.
 \end{equation*} 
 Then $\epsilon_n \longrightarrow 0$, when $n \longrightarrow +\infty$.
\end{lem}

\begin{pf}
We proceed by contradiction. Assume that there exists $\epsilon>0$, a subsequence of $u_n$ (still denoted by $u_n$) and a sequence $t_n \longrightarrow t \in [0,T]$ such that 
  \begin{equation} \label{contra}
 \Big\|    \int_0^{t_n} \big(u_n(\tau) -u(\tau) \big)e^{i(t_n-\tau) H}( K \psi(\tau))d\tau \Big\|_{{\mathcal{H}}^1(\mathbf{R}^3)} \geq \epsilon.
 \end{equation}
Up to a subsequence, we can assume that  for all $n\geq 1$, $t_n \leq t$ or   $t_n \geq t$. We only consider the first case, since the second is similar. By the Minkowski inequality and the unitarity of $e^{i\tau H}$
  \begin{multline*} 
\!\!\! \!\Big\|    \int_0^{t_n}\!\! \!\!\!\!\!\big(u_n(\tau) -u(\tau) \big) \big(   e^{i(t_n-\tau) H} \! - e^{i(t-\tau) H}\big)   ( K \psi(\tau))d\tau \Big\|_{{\mathcal{H}}^1(\mathbf{R}^3)}  \\
   \begin{aligned} 
   &  \leq    \int_0^{t_n} \!\!\!\!\big|u_n(\tau) -u(\tau) \big| \\
   & ~~~~~~~~\times  \Big\|  \big(   e^{i(t_n-\tau) H}- e^{i(t-\tau) H}\big)   ( K \psi(\tau)) \Big\|_{{\mathcal{H}}^1(\mathbf{R}^3)} d\tau\\
      &  =    \int_0^{t_n} \!\!\!\!\big|u_n(\tau) -u(\tau) \big|  \Big\|  \big(   e^{i t_n H}- e^{i t H}\big)   ( K \psi(\tau)) \Big\|_{{\mathcal{H}}^1(\mathbf{R}^3)}d\tau. 
           \end{aligned} 
 \end{multline*}
Then by H\"older
  \begin{eqnarray*} 
\lefteqn{ \Big\|  \!\!\!  \int_0^{t_n} \!\!\!\! \big(u_n(\tau) -u(\tau) \big) \big(   e^{i(t_n-\tau) H}\!\!-\!\! e^{i(t-\tau) H}\!\big)   ( K \psi(\tau))d\tau \Big\|_{{\mathcal{H}}^1(\!\mathbf{R}^3\!)}   }\\
   &  \leq  & \|u_n -u  \|_{L^{1+\epsilon}_T}  \Big\|  \big(   e^{i t_n H}- e^{i t H}\big)   ( K \psi(\tau)) \Big\|_{L^{q_\epsilon}_{\tau \in [0,T]}{\mathcal{H}}^1(\mathbf{R}^3)},           
 \end{eqnarray*}
where $1< q_{\epsilon}< \infty$ is such that $1/(1+\epsilon)+1/q_{\epsilon}=1$. Now, by Lemma~\ref{lem-est}, we have
   \begin{equation}\label{borneL2}
   \| K \psi\|_{L^{q_{\epsilon}}_T{\mathcal{H}}^1(\mathbf{R}^3)}  \leq C_{T,\epsilon} \|\psi \|_{X^1_T} <\infty.
    \end{equation}
 Now we apply   \cite[Lemma 3.2]{ChTh2} (with $d=3$ and $s=1$) together with the previous lines, and we get that
   \begin{equation} \label{small}
 \Big\|    \int_0^{t_n} \!\!\!\!\!\!\big(u_n(\tau) -u(\tau) \big) \big(   e^{i(t_n-\tau) H}- e^{i(t-\tau) H}\big)   ( K \psi(\tau))d\tau \Big\|_{{\mathcal{H}}^1(\mathbf{R}^3)}  
   \end{equation} 
tends to 0 as $n \longrightarrow +\infty$. \medskip

By the Minkowski inequality, the unitarity of $e^{i\tau H}$ and the H\"older inequality
  \begin{multline}  
 \Big\|    \int_{t_n}^t \big(u_n(\tau) -u(\tau) \big)   e^{i(t-\tau) H}   ( K \psi(\tau))d\tau \Big\|_{{\mathcal{H}}^1(\mathbf{R}^3)}  
      \\
          \leq \int_{t_n}^t \big|u_n(\tau) -u(\tau) \big|     \big\|   K \psi(\tau) \big\|_{{\mathcal{H}}^1(\mathbf{R}^3)}d\tau \\
     \begin{aligned}
      &   \leq        \|u_n -u \|_{L^{1+\epsilon}_{T}}  \| K \psi\|_{L^{q_\epsilon}_{\tau \in [t_n,t]}{\mathcal{H}}^1(\mathbf{R}^3)} \\
            &   \leq       |t-t_n|^{1/q_\epsilon}   \|u_n -u \|_{L^{1+\epsilon}_{T}}  \| K \psi\|_{L^{2q_\epsilon}_{T}{\mathcal{H}}^1(\mathbf{R}^3)}. \label{PrT}
            \end{aligned}    
           \end{multline} 
 Using Lemma~\ref{lem-est} and the fact that  $ \|u_n -u \|_{L^{1+\epsilon}_{T}} \leq C$, we deduce that  the term \eqref{PrT}  tends to 0. We combine this with \eqref{small} to deduce 
   \begin{eqnarray}
 \lefteqn{\Big\|   \int_0^{t_n} \big(u_n(\tau) -u(\tau) \big)e^{i(t_n-\tau) H}( K \psi(\tau))d\tau } \nonumber\\
 && - \int_0^t \big(u_n(\tau) -u(\tau) \big)e^{i(t-\tau) H}( K \psi(\tau))d\tau \Big\|_{{\mathcal{H}}^1(\mathbf{R}^3)}  \!\!\!\! \!\!\!\!\!\!\!\!\rightarrow 0.\label{cont2}
 \end{eqnarray}
 
Let us now prove that $\int_0^t \big(u_n(\tau) -u(\tau) \big)e^{i(t-\tau) H}( K \psi(\tau))d\tau$ tends to  0 in ${\mathcal{H}}^1(\mathbf{R}^3)$, to reach a contradiction with  \eqref{contra}. We set $v(\tau)= e^{i(t-\tau) H}( K \psi(\tau))$. Then by the unitarity of $H$, we have $\|v(\tau)\|_{{\mathcal{H}}^1}=\|K \psi(\tau)\|_{{\mathcal{H}}^1}$, thus by \eqref{borneL2}, $v \in L^{q_{\epsilon}}([0,T]; {\mathcal{H}}^1(\mathbf{R}^3))$. We expand $v$ on the Hermite functions $(h_k)_{k\geq 0}$  (which are the eigenfunctions of $H$) which form a Hilbertian basis of~$L^2(\mathbf{R}^3)$  
   \begin{equation*} 
v(\tau,x)=\sum_{k=0}^{+\infty} \alpha_k(\tau) h_k(x),
 \end{equation*}
 so that we have $\displaystyle \|v(\tau, \cdot) \|^2_{{\mathcal{H}}^1}=\sum_{k=0}^{+\infty} (2k+1)|\alpha_k(\tau)|^2$ and 
    \begin{equation*} 
\displaystyle \| v\|_{L_T^{q_\epsilon}{\mathcal{H}}^1}=\bigg[\int_0^T \Big(  \sum_{k=0}^{+\infty} (2k+1)|\alpha_k(\tau)|^2    \big)^{q_{\epsilon}/2}d \tau\bigg]^{1/q_\epsilon}.
 \end{equation*}
 This implies in particular that 
     \begin{equation} \label{consq}
  \Big(  \sum_{k=0}^{+\infty} (2k+1)|\alpha_k(\tau)|^2    \Big)^{q_{\epsilon}/2} \in L^{q_\epsilon}_T.
 \end{equation}
Denote by $\rho= \sup_{n \geq 0} \|u_n -u \|_{L^{1+\epsilon}_{T}} $. We claim that  there exists $M>0$ large enough such that the function $g(\tau,x)=\sum_{k=0}^{M} \alpha_k(\tau) h_k(x)$ satisfies $\| v-g\|_{L^{q_\epsilon}([0,T] ; {\mathcal{H}}^1(\mathbf{R}^3))} \leq \epsilon/(4 \rho)$. Actually, 
    \begin{equation*} 
\displaystyle \| v-g\|_{L_T^{q_\epsilon}{\mathcal{H}}^1}=\bigg[\int_0^T \Big(  \sum_{k=M+1}^{+\infty} (2k+1)|\alpha_k(\tau)|^2    \big)^{q_{\epsilon}/2}d \tau\bigg]^{1/q_\epsilon} 
 \end{equation*}
tends to zero when $M$ tends to $0$, by the Lebesgue theorem and \eqref{consq}, hence the claim.

We have
   \begin{equation*} 
 \int_0^t¸\!\!\!\! \big(u_n(\tau) -u(\tau) \big)g(\tau)d\tau =\sum_{k=0}^{M}   h_k\int_0^t \!\!\!\big(u_n(\tau) -u(\tau) \big) \alpha_k(\tau)d\tau.
 \end{equation*}
Then, by \eqref{consq}, for all $k \geq 0$, $\alpha_k \in  L^{q_\epsilon}_T$, which implies
    \begin{eqnarray*} 
\lefteqn{\big \|\int_0^t \big(u_n(\tau) -u(\tau) \big)g(\tau)d\tau \big\|^2_{{\mathcal{H}}^1(\mathbf{R}^3)}}\\
&&=\sum_{k=0}^{M}  (2k+1)^s\Big|\int_0^t \big(u_n(\tau) -u(\tau) \big) \alpha_k(\tau)d\tau\Big|^2 \longrightarrow 0,
 \end{eqnarray*}
 by the weak convergence of $(u_n)$. Finally, for $n$ large enough,
   \begin{eqnarray*} 
\lefteqn{ \Big\|   \int_0^t \big(u_n(\tau) -u(\tau) \big)v(\tau)d\tau \Big\|_{{\mathcal{H}}^1(\mathbf{R}^3)} }\\
&\leq & \frac{\epsilon}{4\rho}  \big\|u_n -u \big\|_{L^{1+\epsilon}_{T}} +  \Big\|   \int_0^t \big(u_n(\tau) -u(\tau) \big)g(\tau)d\tau \Big\|_{{\mathcal{H}}^1(\mathbf{R}^3)}\\
&\leq &\frac{\epsilon}2
\nonumber \\
 \end{eqnarray*}
which together with \eqref{contra} and \eqref{cont2} gives the contradiction.
\end{pf}

Thanks to \cite[Lemma A.3]{ChTh2} we get
\begin{eqnarray}
\lefteqn{\| \mathcal{N}(\psi,\psi_n)(t) \|_{{\mathcal{H}}^1(\mathbf{R}^3)}} \nonumber\\ 
&\leq & \int_0^t   \|   (\psi-\psi_n)(\psi+\psi_n)  \overline{\psi}  \|_{{\mathcal{H}}^1(\mathbf{R}^3)} d\tau \nonumber\\
&&~~+\int_0^t   \|  (\overline{\psi}-\overline{\psi_n})\psi^2_n  \|_{{\mathcal{H}}^1(\mathbf{R}^3)} d\tau \nonumber \\
&\leq & \int_0^t   \|  z_n  \|_{{\mathcal{H}}^1(\mathbf{R}^3)} \big(  \|  \psi  \|^2_{{\mathcal{W}}^{1,6}} + \|  \psi_n  \|^2_{{\mathcal{W}}^{1,6}}  \big)d\tau. \label{tdl2}
 \end{eqnarray}
 
By \eqref{zn} we have
\begin{multline*}
H^{1/2}z_n (t) =  -i  H^{1/2} \int_0^t \big(u(\tau)-u_n(\tau) \big) e^{i(t-\tau) H}( K \psi)d\tau\\
  -i \int_0^t  u_n(\tau)  e^{i(t-\tau) H}H^{1/2}\big( Kz_n\big)d\tau+ H^{1/2}\mathcal{N}(\psi,\psi_n)(t).
 \end{multline*}
Thus from the Strichartz inequality \eqref{Stri1.1} we deduce
\begin{eqnarray}
 \lefteqn{ \| z_n (t)\|_{{\mathcal{H}}^1(\mathbf{R}^3)}  \leq  \epsilon_n +   \|  u_n H^{1/2}\big( Kz_n\big) \|_{L^2_t L^{6/5}(\mathbf{R}^3)} }\nonumber \\
 &~~~~~~~~~~~~~~~~~~~~~~&+\| \mathcal{N}(\psi,\psi_n)(t) \|_{{\mathcal{H}}^1(\mathbf{R}^3)}.\label{tdl1}
 \end{eqnarray}
 By \eqref{borneKW}
  \begin{equation*}
 \|  H^{1/2}\big( Kz_n\big) \|_{L^{6/5}(\mathbf{R}^3)} = \|  Kz_n \|_{{\mathcal{W}}^{1,6/5}(\mathbf{R}^3)} \leq C  \|  z_n \|_{{\mathcal{H}}^{1}(\mathbf{R}^3)},
 \end{equation*}
 which in turn implies 
  \begin{eqnarray}
 \lefteqn{\|  u_n H^{1/2}\big( Kz_n\big) \|_{L^2_t L^{6/5}(\mathbf{R}^3)} } \nonumber \\
 &~~~~~~~~~~~~~~~~\leq & C \Big(  \int_0^t  |u_n(\tau)|  \|  z_n(\tau) \|^2_{{\mathcal{H}}^{1}(\mathbf{R}^3)} d\tau \Big)^{1/2}.\label{tdl3}
    \end{eqnarray}
  As a conclusion, from \eqref{tdl2}, \eqref{tdl1} and \eqref{tdl3}  we infer
 \begin{multline*} 
   \| z_n (t)\|_{{\mathcal{H}}^1(\mathbf{R}^3)}  \leq  \epsilon_n +  C \Big(  \int_0^t \!\! |u_n(\tau)|  \|  z_n(\tau) \|^2_{{\mathcal{H}}^{1}(\mathbf{R}^3)} d\tau \Big)^{1/2} \!\! +\\
   + \int_0^t   \|  z_n(\tau)  \|_{{\mathcal{H}}^1(\mathbf{R}^3)} \big(  \|  \psi(\tau)  \|^2_{{\mathcal{W}}^{1,6}} + \|  \psi_n(\tau)  \|^2_{{\mathcal{W}}^{1,6}}  \big)d\tau.
   \end{multline*}
  Then by the Gr\"onwall inequality, for all $0 \leq t \leq T$ and~\eqref{boundW16}
   \begin{eqnarray*} 
 \lefteqn{\|z_n(t)\|_{{\mathcal{H}}^1(\mathbf{R}^3)} }\\
 &\leq & \Big( \epsilon_n +  C \Big(  \int_0^t  |u_n(\tau)|  \|  z_n(\tau) \|^2_{{\mathcal{H}}^{1}(\mathbf{R}^3)} d\tau \Big)^{1/2}  \Big) \\
 & &~~ \times  \e^{ c \int_0^t   \big(  \|  \psi(\tau)  \|^2_{{\mathcal{W}}^{1,6}} + \|  \psi_n(\tau)  \|^2_{{\mathcal{W}}^{1,6}}  \big)d\tau} \nonumber\\
 &\leq&C_1(T) \Big( \epsilon_n +   \Big(  \int_0^t  |u_n(\tau)|  \|  z_n(\tau) \|^2_{{\mathcal{H}}^{1}(\mathbf{R}^3)} d\tau \Big)^{1/2}  \Big).  
 \end{eqnarray*} 
We square the previous inequality and get
    \begin{equation*} 
 \|z_n(t)\|^2_{{\mathcal{H}}^1(\mathbf{R}^3)} \leq  2C^2_1(T) \Big( \epsilon^2_n +   \int_0^t \!\! |u_n(\tau)|  \|  z_n(\tau) \|^2_{{\mathcal{H}}^{1}(\mathbf{R}^3)} d\tau  \Big).
 \end{equation*} 
 By the Gr\"onwall inequality again we deduce
 \begin{eqnarray*}
  \lefteqn{\| z_n \|_{L_T^{\infty}{\mathcal{H}}^1(\mathbf{R}^3)} }\\
  & \leq &  2C_1(T) \epsilon_n \exp \big( C_1^2(T)\int_0^T  |u_n(\tau)  |      d\tau  \big) \leq C_2(T) \epsilon_n,
 \end{eqnarray*}
 and this latter quantity tends to 0 when $n \longrightarrow +\infty$.

\section{Conclusion}
This note provides an example of a Ball-Marsden-Slemrod like obstruction to the controllability of a nonlinear partial differential equation with a bilinear control term. The novelty of the result lies  in the unboundedness of the bilinear control term.

The possible relations of this obstruction result and the concepts introduced in \cite{BCC} will be the subject of further investigations in future works. 


 \bibliography{ifacNLS}

\end{document}